\documentclass{article}    
\usepackage{amstext}

\title{Extending the idea of compressed algebra to arbitrary socle-vectors, II: cases of non-existence}  
\author{Fabrizio Zanello\\
Department of Mathematics and Statistics\\
Queen's University\\
Kingston, Ontario\\
K7L 3N6 Canada\\
E-mail: fabrizio@mast.queensu.ca\\
\\}         

\def\qed{$\rlap{$\sqcap$}\sqcup$}
\usepackage{latexsym}

\begin{document}           

{\large

\maketitle                 

{\ }\\

\section{Introduction}

This work is a continuation of the paper $[Za]$, where we introduced {\it generalized compressed algebras} and began their study.\\ Before explaining our new results, let us recall the setting: we work with standard graded artinian algebras, i.e. artinian quotients $A=R/I$ of the polynomial ring $R=k[x_1,...,x_r]$, where the $x_i$'s have degree 1 and $I$ is a homogeneous ideal of $R$. We assume that $k$ is a field of characteristic zero.\\
The {\it h-vector} of $A=\bigoplus_{i=0}^eA_i$ is $h(A)=h=(1,h_1,...,h_e)$, where $h_i=\dim_k A_i$ and $e$ is the last index such that $\dim_k A_e>0$. Since we may suppose that $I$ does not contain non-zero forms of degree 1, $r=h_1$ is defined to be the {\it embedding dimension} ({\it emb.dim.}, in brief) of $A$.\\
The {\it socle} of $A$ is the annihilator of the maximal homogeneous ideal $\overline{m}=(\overline{x_1},...,\overline{x_r})\subset A$, i.e. $soc(A)=\lbrace a\in A {\ } \mid {\ } a\overline{m}=0\rbrace $. Since $soc(A)$ is a homogeneous ideal, we define the {\it socle-vector} of $A$ as $s(A)=s=(s_0,...,s_e)$, where $s_i=\dim_k soc(A)_i$. Notice that $s_0=0$ and $s_e=h_e>0$. The {\it type} of $s$ is type$(s)=\sum_{i=1}^es_i$.\\
The {\it minimum embedding dimension} of a socle-vector $s$ (briefly, \linebreak {\it min.emb.dim.$(s)$}) is defined as the least integer $r$ such that there exists any algebra $A$ with data $(r,s)$. It is easy to see that, if there exists an algebra with data $(r,s)$, then there also exists an algebra with data $(r+1,s)$.\\
We will say that an $h$-vector $h$ is {\it admissible for the pair $(r,s)$} if there exists an algebra $A$ with emb.dim.($A$)$=r$, $s(A)=s$ and $h(A)=h$. When the pair $(r,s)$ is clear from the context, we will simply say that $h$ is {\it admissible}.\\
We place a partial ordering on $h$-vectors of the same length, by defining $h\geq h^{'}$ if every entry of $h$ is greater than or equal to the corresponding entry of $h^{'}$. Using this ordering, given a pair $(r,s)$, two admissible $h$-vectors, $h$ and $h^{'}$, are called {\it comparable} if either $h\geq h^{'}$ or $h^{'}> h$. Otherwise, $h$ and $h^{'}$ are {\it non-comparable}.\\
Finally, given a pair $(r,s)$, define an admissible $h$-vector $h$ as a {\it relative maximum} for the set of all the admissible $h$-vectors if, for every other admissible $h^{'}$ which is comparable with $h$, we have $h\geq h^{'}$.\\
If $h\geq h^{'}$ for {\it every} other admissible $h$-vector $h^{'}$, i.e. if $h$ is the only relative maximum for the set of all the admissible $h$-vectors, we will simply say that $h$ is {\it the maximum}.\\
\\
The problem of finding {\it all} the admissible $h$-vectors for a given pair $(r,s)$ seems very difficult in general, even in the Gorenstein case, i.e. when $s=(0,...,0,1)$. Iarrobino (cf. $[Ia]$) and Fr\"oberg and Laksov (cf. $[FL]$) considered a more restricted question. More precisely, Iarrobino, putting some natural restrictions on a given pair $(r,s)$, showed that any admissible $h$-vector is bounded from above by a certain maximal $h$, and defined an algebra $A$ with the data $(r,s)$ as {\it compressed} if this maximal $h$ satisfies $h=h(A)$; moreover, he proved that, under his hypotheses on $r$ and $s$, there always exists a compressed algebra.\\
This problem of Iarrobino's was taken up again in $[FL]$ by Fr\"oberg and Laksov, who used a different approach.\\
We finally recall the seminal work on compressed algebras, Emsalem and Iarrobino's 1978 article $[EI]$.\\
\\
In our previous paper $[Za]$, we took a more general view and considered the following question: given {\it any} $(r,s)$, is there a maximum $h$ among all the admissible $h$-vectors? If such an $h$ exists, we defined as {\it generalized compressed} any algebra with the data $(r,s,h)$ (see $[Za]$, Def. 2.7). Naturally, this more general definition coincides with Iarrobino's in the cases satisfying his conditions, and, with our generalized definition, we enlarged the set of compressed algebras beyond those found in $[Ia]$ and $[FL]$.\\
Our main contributions in $[Za]$ were: an upper-bound, sharper than that of Fr\"oberg and Laksov, for the admissible $h$-vectors for a given pair $(r,s)$ (cf. $[FL]$, Prop. 4, $i$) and $[Za]$, Thm. A); the theorem that, under certain restrictions on $(r,s)$, our upper-bound is actually achieved by a generalized compressed algebra (cf. $[Za]$, Thm. B). As we saw from some examples, the hypotheses of Theorem B of $[Za]$ cannot (in general) be improved, i.e. under weaker conditions on the pair $(r,s)$, the upper-bound of $[Za]$, Theorem A is not always admissible.\\
\\
Instead, in Section 3 of this paper we will use a result of Cho and Iarrobino on Gorenstein $h$-vectors to exhibit a new class of socle-vectors that admit a generalized compressed algebra (whose $h$-vector is lower than the upper-bound given by $[Za]$, Theorem A). In particular, we will deduce that for every socle-vector $s$ of type 2 there exists a generalized compressed algebra.\\
In Section 4, using Stanley's characterization of Gorenstein $h$-vectors of embedding dimension 3, we will prove the most important result of this paper: there exist pairs $(r,s)$ that do {\it not} admit a generalized compressed algebra, i.e. for these pairs $(r,s)$ the set of all the admissible $h$-vectors has more than one relative maximum.\\
We will also show that (unfortunately??) the scenario can be as bad as we want, even in embedding dimension $r=3$: for every $M>0$ we will construct a pair $(3,s)$ whose set of admissible $h$-vectors has more than $M$ different relative maxima.\\
A question which naturally arises now is the following: what are the pairs $(r,s)$ which admit a generalized compressed algebra?\\
At this stage it seems very difficult to give an answer to this question in the general case. Nevertheless, in Section 5 we will make a first step in this direction, limiting ourselves to a particular class of socle-vectors in embedding dimension 3.\\
\\ 
The results obtained in this paper will be part of the author's Ph.D. dissertation, written at Queen's University (Kingston, Ontario, Canada), under the supervision of Professor A.V. Geramita.
\\
\\
\section{Preliminary results}

Fix $r$ and $s=(s_0=0,s_1,...,s_e)$; from now on we may suppose, to avoid trivial cases, that $r>1$ and $e>1$. Recall that $R=k[x_1,...,x_r]$.\\
\\
{\bf Definition-Remark 2.1.} Following $[FL]$, define, for $d=0,1,...,e$, the integers $$r_d=N(r,d)-N(r,0)s_d-N(r,1)s_{d+1}-...-N(r,e-d)s_e,$$ where $$N(r,d)=\dim_kR_d={r-1+d \choose d}.$$
It is easy to show (cf. $[FL]$) that $r_0<0$, $r_e\geq 0$ and $r_{d+1}>r_d$ for every $d$.\\
Define $b$, then, as the unique index such that $1\leq b\leq e$, $r_b\geq 0$ and $r_{b-1}<0$.\\
\\
Let us now recall briefly the main facts of the theory of Inverse Systems which we will use throughout the paper. For a complete introduction, we refer the reader to $[Ge]$ and $[IK]$.\\
Let $S=k[y_1,...,y_r]$, and consider $S$ as a graded $R$-module where the action of $x_i$ on $S$ is partial differentiation with respect to $y_i$.\\
There is a one-to-one correspondence between artinian algebras $R/I$ and finitely generated $R$-submodules $M$ of $S$, where $I$ is the annihilator of $M$ in $R$ and, conversely, $M$ is the $R$-submodule of $S$ which is annihilated by $I$ (cf. $[Ge]$, Remark 1), p. 17).\\
If $R/I$ has data $(r,s)$, then $M$ is minimally generated by $s_i$ elements of degree $i$, for $i=1,...,e$, and the $h$-vector of $R/I$ is given by the number of linearly independent derivatives in each degree obtained by differentiating the generators of $M$ (cf. $[Ge]$, Remark 2), p. 17).\\
In particular, Gorenstein algebras correspond to cyclic $R$-submodules of $S$.\\
Notice that, given an $R$-submodule $M$ of $S$, if we construct a new $R$-submodule $M^{'}$ by adding $t$ minimal generators in degree $p$ to $M$, then the $h$-vector of the algebra corresponding to $M^{'}$ is unchanged in degrees larger than $p$, and it increases by exactly $t$ in degree $p$.\\
The number $$N(r,d)-r_d=N(r,0)s_d+N(r,1)s_{d+1}+...+N(r,e-d)s_e$$ is an upper-bound for the number of linearly independent derivatives yielded in degree $d$ by the minimal generators of $M$ and, therefore, is also an upper-bound for the $h$-vector of $R/I$. This is the reason for the introduction of the numbers $r_d$.\\
\\
{\bf Proposition 2.2} (Fr\"oberg-Laksov). Let $(r,s)$ be as above, $r\geq $ min.emb.dim.$(s)$. Then an upper-bound for the $h$-vectors admissible for the pair $(r,s)$ is given by $$H=(h_0,h_1,...,h_e),$$ where $$h_i=\min \lbrace N(r,i)-r_i,N(r,i)\rbrace$$ for $i=0,1,...,e$.\\
\\
{\bf Proof.} See $[FL]$, Prop. 4, $i$).{\ }{\ }\qed \\
\\
{\bf Remark 2.3.} A second proof of the proposition follows immediately from our comment about Inverse Systems and the numbers $r_d$. The same upper-bound was already supplied by Iarrobino (cf. $[Ia]$) under the natural restriction $s_1=...=s_{b-1}=0$, where $b$ is as in Definition-Remark 2.1.\\
\\
{\bf Proposition 2.4} (Iarrobino, Fr\"oberg-Laksov). Let $(r,s)$ be as above, $r\geq $ min.emb.dim.$(s)$. If, moreover, $s_1=...=s_{b-1}=0$, then the upper-bound $H$ yielded by Proposition 2.2 is admissible for the pair $(r,s)$.\\
\\
{\bf Proof.} See $[Ia]$, Thm. II A; $[FL]$, Prop. 4, $iv$) and Thm. 14.{\ }{\ }\qed \\
\\
{\bf Definition-Remark 2.5.} Let $n$ and $i$ be positive integers. The {\it i-binomial expansion of n} is $$n_{(i)}={n_i\choose i}+{n_{i-1}\choose i-1}+...+{n_j\choose j},$$ where $n_i>n_{i-1}>...>n_j\geq j\geq 1$.\\
Under these hypotheses, the $i$-binomial expansion of $n$ is unique.\\
Following $[BG]$, define, for any integer $a$, $$(n_{(i)})_{a}^{a}={n_i+a\choose i+a}+{n_{i-1}+a\choose i-1+a}+...+{n_j+a\choose j+a}.$$
\\
A well-known result of Macaulay is:\\
\\
{\bf Theorem 2.6} (Macaulay). Let $h=(h_i)_{i\geq 0}$ be a sequence of non-negative integers, such that $h_0=1$, $h_1=r$ and $h_i=0$ for $i>e$. Then $h$ is the $h$-vector of some standard graded artinian algebra if and only if, for every $1\leq d\leq e-1$, $$h_{d+1}\leq ((h_d)_{(d)})_{+1}^{+1}.$$
\\
{\bf Proof.} See $[St]$.{\ }{\ }\qed \\
\\
{\bf Lemma 2.7} (Bigatti-Geramita). Let $a,b$ be positive integers, $b>1$. Then the smallest integer $s$ such that $a\leq (s_{(b-1)})_{+1}^{+1}$ is $$s=(a_{(b)})_{-1}^{-1}.$$
\\
{\bf Proof.} See $[BG]$, Lemma 3.3.{\ }{\ }\qed \\
\\
{\bf Remark 2.8.} This result yields a lower-bound for the $i$-th entry of an $h$-vector, once the $(i+1)$-st entry is known. In terms of Inverse Systems, it supplies a lower-bound for the number of linearly independent first derivatives of any given set of linearly independent forms of degree $i+1$.\\
\\
An upper-bound, sharper than that of Proposition 2.2, is:\\
\\ 
{\bf Theorem 2.9} ($[Za]$). Let $(r,s)$ be as above, $r\geq $ min.emb.dim.$(s)$. Then an upper-bound $H$ for the $h$-vectors admissible for the pair $(r,s)$ is given by $$H=(h_0,h_1,...,h_e),$$ where $h_0=1$, $h_1=r$ and, inductively, for $2\leq i\leq e$, $$h_i=\min \lbrace ((h_{i-1}-s_{i-1})_{(i-1)})_{+1}^{+1}, N(r,i)-r_i\rbrace .$$
\\
{\bf Proof.} See $[Za]$, Thm. A.{\ }{\ }\qed \\
\\
Comparing the two upper-bounds, we have:\\
\\
{\bf Proposition 2.10} ($[Za]$). The upper-bounds $H$ yielded by Proposition 2.2 and Theorem 2.9 coincide if and only if $s_0=s_1=...=s_{b-2}=0$ and $$s_{b-1}\leq N(r,b-1)-((N(r,b)-r_b)_{(b)})_{-1}^{-1}.$$
Otherwise, Theorem 2.9 yields a sharper $H$.\\
\\
{\bf Proof.} See $[Za]$, Prop. 3.3.{\ }{\ }\qed \\
\\
Before stating the next theorem, notice that we always have $$\max \lbrace N(r,b-1)-(N(r,b)-r_b),0\rbrace \leq N(r,b-1)-((N(r,b)-r_b)_{(b)})_{-1}^{-1},$$
since clearly $((N(r,b)-r_b)_{(b)})_{-1}^{-1}\leq ((N(r,b)-r_b),$ and $r_b>0$, whence $N(r,b-1)-((N(r,b)-r_b)_{(b)})_{-1}^{-1}\geq N(r,b-1)-((N(r,b))_{(b)})_{-1}^{-1}=0$.\\
\\
{\bf Theorem 2.11} ($[Za]$). Fix $(r,s)$ and suppose that $s_0=s_1=...=s_{b-2}=0$. Then, for $$s_{b-1}\leq \max \lbrace N(r,b-1)-(N(r,b)-r_b),0\rbrace ,$$ the upper-bound $H$ of Proposition 2.2 and Theorem 2.9 (which is the same, by Proposition 2.10 and the observation above) is admissible for the pair $(r,s)$.\\
\\
{\bf Proof.} See $[Za]$, Thm. 3.4.{\ }{\ }\qed \\
\\
Recall that a vector $v=(v_0,v_1,...,v_e)$ is called {\it differentiable} if its first difference, $$\Delta v=((\Delta v)_0=1,(\Delta v)_1=v_1-v_0,...,(\Delta v)_e=v_e-v_{e-1}),$$ is an $O$-sequence (i.e., it is the $h$-vector of some standard graded algebra). It is easy to see that if $v$ is differentiable, then $v$ is itself an $O$-sequence.\\
With this definition, we can state two important results about Gorenstein $h$-vectors:\\
\\
{\bf Theorem 2.12} (Stanley). Let $h=(h_0,h_1,...,h_e)$, $h_1=3$. Then $h$ is a Gorenstein $h$-vector if and only if $h$ is symmetric with respect to ${e\over 2}$ and its first half, $(h_0,h_1,...,h_{\lfloor {e\over 2}\rfloor })$, is differentiable.\\
\\
{\bf Proof.} See $[St]$, Thm. 4.2.{\ }{\ }\qed \\
\\
{\bf Theorem 2.13} ($[CI]$). Let $h=(h_0,h_1,...,h_e)$ be an $h$-vector, symmetric with respect to ${e\over 2}$, whose first half, $(h_0,h_1,...,h_{\lfloor {e\over 2}\rfloor })$, is differentiable. Then $h$ is a Gorenstein $h$-vector.\\
\\
{\bf Proof.} See $[CI]$.{\ }{\ }\qed \\
\\
Finally, in Section 5 we will need Theorem B of $[Za]$, which is a generalization of Theorem 2.11 to any socle-vector $s$.\\
\\
{\bf Definition-Remark 2.14.} Fix the pair $(r,s)$, where $r\geq $ min.emb.dim.$(s)$, and let the $h$-vector $H$ be as in Theorem 2.9. Define $c$ as the largest integer such that $h_c$ is generic (i.e. $h_c=N(r,c)$), and $t$ as the largest integer such that $$h_t=((h_{t-1}-s_{t-1})_{(t-1)})_{+1}^{+1}<N(r,t)-r_t,$$ where we set $(1_{(0)})_{+1}^{+1}=r$ and $((h_{-1}-s_{-1})_{(-1)})_{+1}^{+1}=1$, in order to avoid pathological cases.\\
It is easy to see that $0\leq t\leq e-1$ and $1\leq c\leq t+1$.\\
\\
{\bf Theorem 2.15} ($[Za]$). Let $(r,s)$ be as above, $r\geq $ min.emb.dim.$(s)$, and the upper-bound $H$ given by Theorem 2.9. Then $H$ is admissible in the following cases:
\begin{itemize}
\item[  $i$)] $c=t+1$;
\item[  $ii$)] $c=t$ and $s_c\leq \max \lbrace N(r,c)-h_{c+1},0\rbrace $;
\item[  $iii$)] $c\leq t-1$ and $s_c\geq N(r,c)-c$.
\end{itemize}
{\ }\\
{\bf Proof.} See $[Za]$, Thm. B.{\ }{\ }\qed \\
\\
\\
\section{Generalized compressed algebras: a special class of socle-vectors}

In this section we study the existence of generalized compressed algebras having embedding dimension $r$ for the class of socle-vectors $s=(s_0=0,s_1,...,s_e)$, where $s_e=1$ and only one more entry, say $s_p$, is non-zero, making use of Inverse Systems and the properties of Gorenstein $h$-vectors.\\
In our Theorem 3.1 below the reader will notice the restriction $1\leq s_p\leq r-1$. That this condition is necessary for the conclusion of the theorem will be shown in Example 4.1. In fact, we will see in the next section that, already in emb.dim. 3, the case $s_p\geq r$ leads to a completely different scenario. Also in emb.dim. $r\geq 4$, as we will see from Example 3.4, the study of the case $s_p\geq r$ would be very different, requiring knowledge of Gorenstein $h$-vectors beyond the state of the art.\\
\\
Notice that, under the hypotheses on $s$ that we made above, i.e. that the socle is concentrated in degrees $p$ and $e$ and $s_e=1$, by Definition 2.1 we have $$r_i=N(r,i)-N(r,p-i)s_p-N(r,e-i)$$ for $0\leq i\leq p$, and $$r_i=N(r,i)-N(r,e-i)$$ for $p<i\leq e.$\\
\\
{\bf Theorem 3.1.} Fix $(r,s=(s_0=0,...,s_p,...,s_e))$, $r\geq $ min.emb.dim.$(s)$, where $1\leq s_p\leq r-1$ for some $p<e$, $s_e=1$ and $s_i=0$ otherwise. Then: 1). There exists a generalized compressed algebra for the pair $(r,s)$.\\
Moreover, if we let $H=(1,h_1,...,h_e)$ denote the $h$-vector of this generalized compressed algebra, then:\\
2). a). If $p\geq {e\over 2}$, then $H$ is the upper-bound of Proposition 2.2, i.e. $$h_i=\min \lbrace N(r,i)-r_i,N(r,i)\rbrace $$ for $i=0,1,...,e$;\\
b). If $p<{e\over 2}$, then $H$ is as follows:
\[h_i=N(r,i) \text{ for } i=0,...,p,\] 
\[h_{p+a}=((N(r,p)-s_p)_{(p)})_a^a \text{ for } a=1,...,\lfloor {e\over 2}\rfloor -p,\]
\[h_i=h_{e-i} \text{ for } i>\lfloor {e\over 2}\rfloor , i\neq e-p\]
and
\[h_{e-p}=h_p-s_p=N(r,p)-s_p.\]
\\
Since in the {\it level} case (i.e. when the socle is concentrated only in the last degree) the existence of (generalized) compressed algebras is well-known (e.g., it follows from Proposition 2.4), it is immediate from Theorem 3.1 that:\\
\\
{\bf Corollary 3.2.} Fix $(r,s)$, $r\geq $ min.emb.dim.$(s)$. If type$(s)=2$, then there always exists a generalized compressed algebra (and its $h$-vector is given by Theorem 3.1).\\
\\
Since in emb.dim. 2 we already know that a generalized compressed algebra always exists (see $[Za]$, Prop. 3.12), it also follows from Theorem 3.1 that:\\
\\
{\bf Corollary 3.3.} Fix $(r,s)$, $r\geq $ min.emb.dim.$(s)$. If type$(s)=3$, with $s_e=1$ and $s_p=2$ for some $p<e$, then there exists a generalized compressed algebra (and its $h$-vector is given by Theorem 3.1).\\
\\
{\bf Proof of Theorem 3.1.} {\it Case $p>{e\over 2}$}. We first determine the integer $b$ associated to $(r,s)$ (see Definition-Remark 2.1). Consider $r_p$; as we have observed before the statement of the theorem, $$r_p=N(r,p)-s_p-N(r,e-p).$$
Since $p>e-p$, then $N(r,p)\geq N(r,e-p+1)$. Furthermore, since $s_p\leq r-1$, we obtain $$r_p\geq N(r,e-p+1)-N(r,e-p)-(r-1)=N(r-1,e-p+1)-(r-1),$$ which is greater than or equal to 0, since $e-p+1\geq 1.$ Therefore $b\leq p$. Thus $s_1=...=s_{b-1}=0$, whence it follows, by Proposition 2.4, that the upper-bound $H$ of Proposition 2.2 is admissible.\\
{\it Case $p={e\over 2}$}. Of course now $e$ is even. As in the previous case, we seek the integer $b$ associated to $(r,s)$. As in the observation before the statement, $$r_p=r_{{e\over 2}}=N(r,{e\over 2})-s_p-N(r,{e\over 2})=-s_p<0$$ and $$r_{p+1}=r_{{e\over 2}+1}=N(r,{e\over 2}+1)-N(r,{e\over 2}-1)>0.$$
Hence $b={e\over 2}+1$, i.e. $p=b-1$, and therefore $s_1=...=s_{b-2}=0$.\\
We now want to show that $s_{b-1}\leq \max \lbrace 0,N(r,b-1)-(N(r,b)-r_b)\rbrace ,$ in order to use Theorem 2.11. $$\max \lbrace 0,N(r,b-1)-(N(r,b)-r_b)\rbrace $$$$=\max \lbrace 0,N(r,{e\over 2})-N(r,{e\over 2}+1)+N(r,{e\over 2}+1)-N(r,{e\over 2}-1)\rbrace $$$$=\max \lbrace 0,N(r,{e\over 2})-N(r,{e\over 2}-1)\rbrace =N(r,{e\over 2})-N(r,{e\over 2}-1)$$$$=N(r-1,{e\over 2})\geq r-1\geq s_p=s_{b-1}.$$ Therefore we are in the hypotheses of Theorem 2.11, and it follows that the upper-bound $H$ of Proposition 2.2 is admissible.\\
{\it Case $p<{e\over 2}$}. By Inverse Systems, any $h$-vector $h$ which is admissible for our pair $(r,s)$ must have the entries of degrees $p+1,...,e$ equal to the number of linearly independent derivatives of some form $F$ of degree $e$. Since the $h$-vector of $R/Ann(F)$ is Gorenstein and therefore symmetric, in degrees $p$ and $e-p$ it must be less than or equal to $N(r,p)-s_p$, in order to leave room for $s_p$ more linearly independent forms of degree $p$; thus, by the symmetry of the Gorenstein $h$-vectors and Theorem 2.6, it follows that the $H$ described in the statement for $p<{e\over 2}$ is an upper-bound for the admissible $h$-vectors for our pair $(r,s)$.\\
It remains to show that $H$ is admissible. Notice that this is true if the symmetric $h$-vector whose first half is given by $$v_0=N(r,0), v_1=N(r,1), ...,v_{p-1}=N(r,p-1),$$ $$v_p=N(r,p)-s_p, v_{p+1}=((N(r,p)-s_p)_{(p)})_1^1,$$ $$v_{p+2}=((N(r,p)-s_p)_{(p)})_2^2, ...,v_{\lfloor {e\over 2}\rfloor }=((N(r,p)-s_p)_{(p)})_{\lfloor {e\over 2}\rfloor -p}^{\lfloor {e\over 2}\rfloor -p}$$ is a Gorenstein $h$-vector. To see why this is so, observe that, if this is a Gorenstein $h$-vector, then there exists a form of degree $e$ which yields $N(r,p)-s_p$ derivatives in degree $p$; therefore, if we now add $s_p$ generators in that degree to our Inverse System, we immediately obtain $H$.\\
We turn now to the proof that the vector $$v=(v_0,v_1,...,v_{\lfloor {e\over 2}\rfloor })$$ is the first half of a Gorenstein $h$-vector. Since, by Theorem 2.13, any symmetric $h$-vector whose first half is differentiable is a Gorenstein $h$-vector, it is enough to show that the vector $v$ is differentiable.\\
Let $$\Delta v=((\Delta v)_0=1,(\Delta v)_1=v_1-v_0,...,(\Delta v)_{\lfloor {e\over 2}\rfloor }=v_{\lfloor {e\over 2}\rfloor }-v_{\lfloor {e\over 2}\rfloor -1})$$ be the first difference of $v$. Then we have:
\begin{itemize}
\item[ $i)$] $$(\Delta v)_i=N(r,i)-N(r,i-1)=N(r-1,i)$$ for $i=0,...,p-1$;
\item[ $ii)$] $$(\Delta v)_p=N(r,p)-s_p-N(r,p-1)=N(r-1,p)-s_p;$$
\item[ $iii)$] $$(\Delta v)_{p+a}=((N(r,p)-s_p)_{(p)})_a^a-((N(r,p)-s_p)_{(p)})_{a-1}^{a-1}$$ for $a=1,...,\lfloor {e\over 2}\rfloor -p$.
\end{itemize}
From $i)$ and $ii)$ we immediately obtain that $\Delta v$ is an $O$-sequence up to degree $p$. Thus, by $ii)$ and $iii)$, it is enough to show that \begin{equation}\label{p}((N(r-1,p)-s_p)_{(p)})_1^1\geq ((N(r,p)-s_p)_{(p)})_1^1-(N(r,p)-s_p)\end{equation} and that, for $j\geq 1$,\\
\begin{equation}\label{j}\matrix{(((N(r,p)-s_p)_{(p)})_j^j-((N(r,p)-s_p)_{(p)})_{j-1}^{j-1})_1^1\cr \cr
\geq ((N(r,p)-s_p)_{(p)})_{j+1}^{j+1}-((N(r,p)-s_p)_{(p)})_j^j.\cr}\end{equation}
{\ }\\ We will, in fact, show that equality holds in (\ref{p}) and (\ref{j}).\\
Since $$N(r,p)={r+p-1\choose p}={r+p-2\choose p}+{r+p-3\choose p-1}$$ $$+...+{r\choose 2}+{r-1\choose 1}+1$$ and $1\leq s_p\leq r-1$, we have that \begin{equation}\label{r}(N(r,p)-s_p)_{(p)}={r+p-2\choose p}+...+{r\choose 2}+{r-s_p\choose 1}\end{equation} and, in a similar fashion,\\
\begin{equation}\label{r-1}\matrix{(N(r-1,p)-s_p)_{(p)}\cr \cr
={r+p-3\choose p}+...+{r-1\choose 2}+{r-1-s_p\choose 1}.\cr}\end{equation}
{\ }\\ (Note that, if $s_p=r-1$, then the r.h.s. of (\ref{r-1}) finishes with ${r-1\choose 2}$. So, to continue the argument, we should really separate the two cases $s_p=r-1$ and $s_p<r-1$. We will continue the proof only for the case $s_p<r-1$, the other one having no different features).\\
Therefore, by (\ref{r}) and (\ref{r-1}), showing equality in (\ref{p}) is equivalent to showing that\\
\begin{equation}\label{bi}\matrix{{r+p-2\choose p+1}+{r+p-3\choose p}+...+{r\choose 3}+{r-s_p\choose 2}\cr \cr
={r+p-1\choose p+1}+{r+p-2\choose p}+...+{r+1\choose 3}+{r+1-s_p\choose 2}-{r-1+p\choose p}+s_p.\cr}\end{equation}
{\ }\\ It is easy  to check that, for any integers $a$ and $b$ such that $a-1\geq b\geq 1$, we have the equality \begin{equation}\label{fab}{a\choose b}-{a-1\choose b}={a-1\choose b-1}.\end{equation}
By (\ref{fab}), a simple calculation shows that (\ref{bi}) can be rewritten as $${r-1+p\choose p}-s_p={r+p-2\choose p}+{r+p-3\choose p-1}+...+{r\choose 2}+{r-s_p\choose 1},$$ which is true by (\ref{r}). This proves equality in (\ref{p}).\\
Now we want to show equality in (\ref{j}). By (\ref{r}) and (\ref{fab}), we can see that the l.h.s. of (\ref{j}) is equal to \begin{equation}\label{hh}(({r+p+j-2\choose p+j}-{r+j-1\choose j+1}+{r+j-s_p-1\choose j+1})_{(p+j)})_1^1.\end{equation}
It is easy to check that $${r+p+j-2\choose p+j}-{r+j-1\choose j+1}$$$$={r+p+j-3\choose p+j}+{r+p+j-4\choose p+j-1}+...+{r+j\choose j+3}+{r+j-1\choose j+2}.$$
Therefore (\ref{hh}) is easily seen to equal
\begin{equation}\label{casa}{r+p+j-1\choose p+j+1}-{r+j\choose j+2}+{r-s_p+j\choose j+2}.\end{equation}
Similarly, the r.h.s. of (\ref{j}), by (\ref{r}) and (\ref{fab}), becomes $${r+p+j-1\choose p+j+1}-{r+j\choose j+2}+{r-s_p+j\choose j+2},$$
which is equal to (\ref{casa}). This completes the proof of the theorem.{\ }{\ }\qed \\
\\
{\bf Example 3.4.} Now we will see that Cho and Iarrobino's theorem (see Theorem 2.13), which was a fundamental tool in our proof of Theorem 3.1, is no longer useful if we drop the hypothesis $s_p\leq r-1$.\\
For instance, let $r=4$, $s=(0,0,0,4,0,0,0,0,1)$. Then, reasoning as in Theorem 3.1, an upper-bound for the admissible $h$-vectors for the pair $(r,s)$ is $$H=(1,4,10,20,25,16,10,4,1).$$ But the first difference of $(1,4,10,16,25)$ is $(1,3,6,6,9)$, which is not an $O$-sequence; therefore Theorem 2.13 is not applicable. (Actually, it can be shown that the $h$-vector $H$ is not admissible. However, we still do not know if there exists a generalized compressed algebra for this pair $(r,s)$).\\
The problem in studying the case $s_p\geq r$ for $r\geq 4$ is that, today, very little is known about Gorenstein $h$-vectors whose first half is not differentiable.\\
In emb.dim. $r=3$, instead, we will see in the next section how bad the scenario can be for $s_p\geq r$.\\
\\
\\
\section{Generalized compressed algebras: cases of non-existence}

In this section, as we already mentioned, we will show that there are cases where a generalized compressed algebra does not exist, and that, moreover, the way this pathology occurs can be arbitrarily bad. To do this, we will make strong use of both Inverse Systems and Stanley's characterization of Gorenstein $h$-vectors of emb.dim. 3 (see Theorem 2.12).\\
\\
{\bf Example 4.1.} Let $r=3$, $s=(0,0,0,3,0,0,0,0,1)$. Then, for any $h$-vector $h=(h_0,h_1,...,h_8)$ admissible for this pair $(r,s)$, by Inverse Systems and Macaulay's theorem (2.6), we must have $h_5\leq 7$; in fact, for $i\geq 4$, the $h_i$'s are given by the final part of a Gorenstein (hence symmetric) $h$-vector whose entry of degree 3 must be less than or equal to 7, in order to allow $s_3=3$.\\
A Gorenstein $h$-vector starting with $(1,3,6,7)$ has first difference starting with $(1,2,3,1)$. Then, by Theorem 2.12, since $(1_{(3)})_1^1=1$, the maximal entry of this $h$-vector in degree 4 can be $7+1=8$, and $(1,3,6,7,8,7,6,3,1)$ is a Gorenstein $h$-vector. Using Inverse Systems, we know that there is a form of degree 8 which yields this $h$-vector. Since $s_3=3$, if now we add 3 generators in degree 3 to our Inverse System, we immediately obtain $H^{'}=(1,3,6,10,8,7,6,3,1)$, which therefore is admissible and is also a relative maximum for the admissible $h$-vectors for our pair $(r,s)$.\\
Similarly, by Stanley's theorem, $(1,3,5,7)$ has maximal growth equal to $7+2=9$ in degree 4, and $(1,3,5,7,9,7,5,3,1)$ is a Gorenstein $h$-vector. Therefore, reasoning as above, we can see that also $H^{''}=(1,3,6,10,9,7,5,3,1)$ is admissible and is a relative maximum for the admissible $h$-vectors for our pair $(r,s)$.\\
Since $H^{'}$ and $H^{''}$ are different, it follows that the pair $(r,s)$ of this example admits no generalized compressed algebra.\\
It also can be checked that $H^{'}$ and $H^{''}$ are the only two relative maxima for the admissible $h$-vectors for this pair $(r,s)$. We will see from the next theorem that things may be even worse: in fact the number of different relative maximal $h$-vectors for any pair $(r,s)$ is not bounded from above, even in emb.dim. 3.\\
\\
{\bf Theorem 4.2.} For every $M>0$ there exists a pair $(3,s)$ such that the set of its admissible $h$-vectors has more than $M$ different relative maxima.\\
\\
{\bf Proof.} Fix $n\geq 2$. Let $r=3$ and $s=(s_0,s_1,...,s_e)$, where $e=4n$, $s_{4n}=1$, $s_{2n-1}=2n+1$ and $s_i=0$ otherwise. Consider the $n$ vectors of the form $$H=(1,3,6,...,N(3,2n-1),h_{2n},h_{2n+1},h_{2n+2},N(3,2n-3),...,6,3,1),$$ where the tuple $(h_{2n},h_{2n+1},h_{2n+2})$ assumes the values: $$(N(3,2n-2)-1, N(3,2n-2)-1,N(3,2n-2)-1),$$ $$(N(3,2n-2),N(3,2n-2)-1,N(3,2n-2)-2),$$ $$(N(3,2n-2)+1,N(3,2n-2)-1,N(3,2n-2)-3),$$ $$\vdots $$ $$(N(3,2n-2)+n-2,N(3,2n-2)-1,N(3,2n-2)-n).$$ These are easily seen to be $h$-vectors of artinian algebras, by Macaulay's theorem, since they are generic up to degree $2n-1$ and then they are non-increasing.\\
We want to show that all these $n$ $h$-vectors are relative maxima in the set of the admissible $h$-vectors for the pair $(3,s)$ defined above.\\
First we show that they themselves are admissible.\\
Notice that the midpoint of these $h$-vectors occurs in degree $2n$, and that \begin{equation}\label{gor}h_{2n+1}+s_{2n-1}=N(3,2n-2)-1+(2n+1)=N(3,2n-1),\end{equation} which is equal to the (generic) entry of degree $2n-1$ of all of our $n$ $h$-vectors.\\
Hence, showing that these $h$-vectors are admissible means showing that the $n$ symmetric $h$-vectors whose first half is given by $$v=(1,3,6,...,N(3,2n-3),h_{2n+2},h_{2n+1},h_{2n})$$ are Gorenstein $h$-vectors.\\
In fact, using Inverse Systems, each of these Gorenstein $h$-vectors is generated by a form of degree $4n$. Hence, if now we add $s_{2n-1}=2n+1$ generators in degree $2n-1$ to our Inverse Systems, by (\ref{gor}) we easily obtain the $n$ $h$-vectors above.\\
Therefore, by Theorem 2.12, it is enough to show that $v$ is differentiable for each of the $n$ choices of $(h_{2n},h_{2n+1},h_{2n+2})$ described above. The first difference of $v$ is $$\Delta (v)=(\Delta (v)_0=1,\Delta (v)_1=2,...,\Delta (v)_{2n-3}=2n-2,$$$$\Delta (v)_{2n-2}=h_{2n+2}-N(3,2n-3),\Delta (v)_{2n-1}=h_{2n+1}-h_{2n+2},$$$$\Delta (v)_{2n}=h_{2n}-h_{2n+1}).$$
Notice that $\Delta (v)$ is generic up to degree $2n-3$. In order to prove that it is an $O$-sequence, it suffices to show that it does not increase from degree $2n-3$ to degree $2n$. We have:
$$\Delta (v)_{2n-2}=h_{2n+2}-N(3,2n-3)\leq (N(3,2n-2)-1)-N(3,2n-3)$$$$=2n-2=\Delta (v)_{2n-3}.$$
$\Delta (v)_{2n-2}\geq \Delta (v)_{2n-1}$ is clearly equivalent to $2h_{2n+2}\geq N(3,2n-2)-1+N(3,2n-3)$, and this is easy to check, since $h_{2n+2}\geq N(3,2n-2)-n$.\\
Finally, by definition, we immediately have $\Delta (v)_{2n-1}=\Delta (v)_{2n}$.
Hence we have shown that $\Delta (v)$ is an $O$-sequence, and therefore that the $n$ $h$-vectors $H$ defined above are admissible for our pair $(r,s)$.\\
It remains to show that they are relative maxima. The entries up to $N(3,2n-1)$ are naturally maximal, and the same holds for the entries from $N(3,2n-3)$ on, since, by Proposition 2.2, they are the final entries of an upper-bound for the Gorenstein $h$-vectors.\\
The entry $h_{2n+1}=N(3,2n-2)-1$ is also maximal. In fact, using Inverse Systems, the second half of our $h$-vectors $H$ is generated by a form of degree $4n$, which, by symmetry, has $h_{2n+1}$ derivatives in degree $2n-1$. Since, in our Inverse Systems for $H$, we also have $s_{2n-1}=2n+1$ generators of degree $2n-1$, it immediately follows that $$h_{2n+1}\leq N(3,2n-1)-(2n+1)=N(3,2n-2)-1;$$ hence the entry $h_{2n+1}=N(3,2n-2)-1$ is maximal.\\
Now consider the first tuple $(h_{2n},h_{2n+1},h_{2n+2}),$ i.e. $$(N(3,2n-2)-1, N(3,2n-2)-1,N(3,2n-2)-1).$$
Using Inverse Systems as above, since, by Theorem 2.12, Gorenstein $h$-vectors in emb.dim. 3 are unimodal (i.e. they do not increase once they start decreasing), we have at once that the value $N(3,2n-2)-1$ is maximal for $h_{2n+2}$. Hence, by Stanley's theorem, we obtain $h_{2n}=N(3,2n-2)-1$. Therefore it follows that the $h$-vector $H$ corresponding to the first tuple is a relative maximum.\\
Now, arguing as above, by Inverse Systems and Theorem 2.12 it is easy to see that also the remaining $n-1$ $h$-vectors $H$ described above are relative maxima.\\
Therefore, we have shown that, for each $n\geq 2$, the set of the admissible $h$-vectors for the pair $(3,s)$ defined above has at least $n$ different relative maxima (actually, it can be seen that there are more than $n$ relative maxima for $n\geq 4$). This clearly completes the proof of the theorem.{\ }{\ }\qed \\
\\
\\
\section{Generalized compressed algebras: characterization in a particular case}

As we already stated in the Introduction, it would be very interesting to have a complete characterization of the pairs $(r,s)$ which admit a generalized compressed algebra. In general, today, this problem seems very hard to attack. However, in embedding dimension $r=3$, there are cases where we can be successful.\\
\\
{\bf Theorem 5.1.} Let $r=3$, $s=(0,...,s_p,...,s_e)$, $3\geq $ min.emb.dim.$(s)$. Suppose that $s_p\geq 1$ for some $p<e$, $s_e=1$ and $s_i=0$ otherwise. Then there exists a generalized compressed algebra for the pair $(3,s)$ if and only if one of the following holds:
\begin{itemize}
\item[ $i)$]  $p<\lfloor {e\over 2}\rfloor $ and either $s_p\leq 2$ or $s_p\geq N(3,p)-p$;
\item[ $ii)$] $p\geq \lfloor {e\over 2}\rfloor $.
\end{itemize}
Equivalently, there does not exist a generalized compressed algebra for the pair $(3,s)$ above if and only if: $$p<\lfloor {e\over 2}\rfloor \mbox{{\ }and {\ }} 3\leq s_p<N(3,p)-p.$$
\\
{\bf Proof.} We begin by showing that there exists a generalized compressed algebra in all the cases mentioned in the statement.\\
$i)$. Let $p<\lfloor {e\over 2}\rfloor $. If $s_p\leq 2$, then there exists a generalized compressed algebra by Theorem 3.1, 2) b).\\
Now consider $s_p\geq N(3,p)-p$. Since $N(3,p)-s_p\leq p$, we have $((N(3,p)-s_p)_{(j)})_1^1=N(3,p)-s_p$ for every $j\geq p$. Therefore the upper-bound $H$ given by Theorem 2.9 for this pair $(3,s)$ is $$(1,h_1,...,h_e),$$ where 
$$
h_i=\cases{N(3,i),&if $1\leq i\leq p$\cr
	\min \lbrace N(3,p)-s_p,N(3,e-i)\rbrace ,&if $p+1\leq i\leq e$.\cr}
$$
In particular, we immediately have that $h_i=N(3,p)-s_p$ for $p+1\leq i\leq e-p$.\\
With notation as in Definition-Remark 2.14, it is easy to see that $c=p$ and $t\geq \lfloor {e\over 2}\rfloor $. Hence $c\leq t-1$; thus, by Theorem 2.15, $iii$), there exists a generalized compressed algebra for the pair $(3,s)$.\\
$ii)$. Let $p\geq \lfloor {e\over 2}\rfloor $. It follows that $p+1>e-p-1$, and therefore (see Definition-Remark 2.1),\\ $r_{p+1}=N(3,p+1)-N(3,e-p-1)={2+p+1\choose 2}-{2+e-p-1\choose 2}>0$. Thus $b\leq p+1$, i.e. $p\geq b-1$.\\
We have that $r_p=N(3,p)-s_p-N(3,e-p)$. If $r_p\geq 0$, then $p\geq b$ and the existence of a generalized compressed algebra follows from Proposition 2.4.\\
Thus, suppose $r_p<0$, i.e. $s_p>N(3,p)-N(3,e-p)$.\\
{\it Claim}. The maximal $h$-vector for our pair $(3,s)$ is $$H=(1,h_1,...,h_e),$$ where
$$
h_i=\cases{N(3,i),&if $1\leq i\leq p$\cr
	\min \lbrace N(3,i),N(3,p)-s_p\rbrace ,&if $p+1\leq i\leq e$.\cr}
$$
{\it Proof of claim.} We first show that $H$ is an upper-bound. By Inverse Systems, the entries of degree $i$, for $i\geq p+1$, of any $h$-vector admissible for this pair $(3,s)$ are given by the number of linearly independent derivatives (of order $e-i$) of a form $F$ of degree $e$; moreover, $F$ must have at most $N(3,p)-s_p$ derivatives of order $e-p$, in order to leave room for $s_p$ more linearly independent forms of degree $p$. Since $p\geq \lfloor {e\over 2}\rfloor $ and Gorenstein $h$-vectors of embedding dimension 3 are unimodal, we must have $h_i\leq N(3,p)-s_p$ for every $i\geq p+1$. This immediately shows that $H$ is an upper-bound for the admissible $h$-vectors for the pair $(3,s)$.\\
Now we want to show that $H$ is admissible. Reasoning as above, it is enough to show that the entries of degrees $i\geq p+1$ are the final entries of a Gorenstein $h$-vector, but this is clear by construction. Hence the proof of the claim is complete.\\
Therefore we have shown $ii)$.\\
Let us assume, for the rest of the proof, that $p<\lfloor {e\over 2}\rfloor $ and $3\leq s_p<N(3,p)-p$. It remains to show that, under these hypotheses, there exists no generalized compressed algebra for the pair $(3,s)$.\\
We break up the range $[3,N(3,p)-p]$ into $[3,p+1]$ and $[p+1,N(3,p)-p]$ (since clearly $p+1\leq N(3,p)-p$).\\
Define $a=p+1-s_p.$ Notice that $a\leq p-2$.\\
Let us consider first the case $3\leq s_p\leq p+1$, i.e. $a\geq 0$.\\
{\it Claim.} A relative maximum for this pair $(3,s)$ is $$H=(1,h_1,...,h_e),$$ where 
\[h_i=N(3,e-i) \text{ for } i\geq e-p+1,\] 
\[h_{e-p-j}=N(3,p)-s_p+aj \text{ for } j=0,1,...,\lfloor {e\over 2}\rfloor -p,\] 
\[h_i=h_{e-i} \text{ for } i\leq \lfloor {e\over 2}\rfloor , i\neq p\]
and
\[h_p=N(3,p)=h_{e-p}+s_p.\]
{\it Proof of claim.} We will just sketch the argument, since it closely follows that for Theorem 4.2.\\
We first show that $H$ is admissible. By Inverse Systems, it is enough to show that the second half of $H$, written the other way around, is the first half, say $v$, of a Gorenstein $h$-vector, and this easily follows from Stanley's theorem, since as soon as $v$ is no longer generic, its first difference, $\Delta (v)$, by construction, becomes constantly equal to $a$.\\
To see why $H$ is also a relative maximum, just observe that the entries of degree $i\leq p$ and $i\geq e-p$ must be maximal, and that $\Delta (v)=(1,2,...,p,a,a,...,a,0)$. Since $a\leq p-2$, it easily follows that $H$ is a relative maximum, and the proof of the claim is complete.\\
In order to show that the pair $(3,s)$ above, for $3\leq s_p\leq p+1$, admits no generalized compressed algebra, it is enough to exhibit an admissible $h$-vector for the pair $(3,s)$ which is not comparable with $H$.\\
{\it Claim.} Such an $h$-vector is $$H^{'}=(1,h_1^{'},...,h_e^{'}),$$ where 
\[h_i^{'}=h_i \text{ for } i\leq p \text{ and } i\geq e-p+2,\] 
\[h_{e-p+1}^{'}=h_{e-p+1}-1=N(3,p-1)-1,\text{\ }\] 
\[h_{e-p-j}^{'}=N(3,p)-s_p+(a+1)j \text{ for } j=0,1,...,\lfloor {e\over 2}\rfloor -p,\]
\[h_i^{'}=h_{e-i}^{'} \text{ for } i\leq \lfloor {e\over 2}\rfloor , i\neq p\]
and
\[h_p^{'}=N(3,p)=h_{e-p}^{'}+s_p.\text{\ } \]
{\it Proof of claim.} Notice that $H$ and $H^{'}$ are non-comparable, since $h_{e-p+1}^{'}<h_{e-p+1}$ and $h_{e-p-j}^{'}>h_{e-p-j}$ for $j=1,...,\lfloor {e\over 2}\rfloor -p$. Therefore it remains to show that $H^{'}$ is admissible.\\
By the same argument we used above, it is enough to show that the second half of $H^{'}$, written the other way around, is the first half, $v$, of a Gorenstein $h$-vector, i.e., by Stanley's theorem, that the first difference, $\Delta (v)$, of $v$ is an $O$-sequence.\\
We have $$\Delta (v)=(1,2,...,p-1,p-1,a+1,a+1,...,a+1,0).$$ Hence it suffices to show that $\Delta (v)$ is non-increasing after degree $p-1$, i.e. that $p-1\geq a+1$. But, since $s_p\geq 3$, we immediately have $a+1=p+1-s_p\leq p-1$, and therefore the proof of the claim is complete.\\
This finishes the case $3\leq s_p\leq p+1$, i.e. $a\geq 0$.\\
It remains to show that also for $p<\lfloor {e\over 2}\rfloor $ and $p+1\leq s_p<N(3,p)-p$ (i.e. $a\leq 0$) there exists no generalized compressed algebra for the pair $(3,s)$. Again, it is enough to exhibit a relative maximum $H$ and an admissible $h$-vector $H^{'}$ which is not comparable with $H$.\\
{\it Claim.} A relative maximum for our pair $(3,s)$ is $$H=(1,h_1,...,h_e),$$ where 
$$
h_i=\cases{N(3,i),&if $1\leq i\leq p$\cr
	\min \lbrace N(3,p)-s_p,N(3,e-i)\rbrace ,&if $p+1\leq i\leq e$.\cr}
$$
{\it Proof of claim}. This proof is similar to the previous one and will be omitted.\\
It remains to find an admissible $h$-vector $H^{'}$ which is not comparable with $H$.\\
{\it Claim.} Such an $h$-vector is $$H^{'}=(1,h_1^{'},...,h_e^{'}),$$ where 
\[
h_i^{'}=N(3,i) \text{ for } 1\leq i\leq p,
\]
\[
h_{e-p+j}^{'}=\min \lbrace N(3,p)-s_p-j,N(3,p-j)\rbrace \text{ for }j=0,1,...,p,
\]
\[ h_{e-p-j}^{'}=N(3,p)-s_p+j\text{ for }j=0,1,...,\lfloor {e\over 2}\rfloor -p
\] and
\[
h_i^{'}=h_{e-i}^{'}\text{ for }i=p+1,...,\lfloor {e\over 2}\rfloor .
\]
{\it Proof of claim.} First of all, notice that $H$ and $H^{'}$ are non-comparable. In fact it is easy to see that, since $p<\lfloor {e\over 2}\rfloor $, we have $h_{e-\lfloor {e\over 2}\rfloor }^{'}>h_{e-\lfloor {e\over 2}\rfloor }$, and $H$ is a relative maximum for the pair $(3,s)$.\\
Hence, it remains to show that $H^{'}$ is admissible. As usual, it suffices to prove that the second half of $H^{'}$, written the other way around, is equal to the first half, $v$, of a Gorenstein $h$-vector of embedding dimension 3, i.e. that the first difference, $\Delta (v)$, of $v$ is an $O$-sequence.\\
Let $q$ be the smallest integer such that $$N(3,p)-s_p-q\geq N(3,p-q).$$ It is easy to see that, since $p+1\leq s_p\leq N(3,p)-p$, $q$ is uniquely determined in the range $[2,p]$.\\
We have $$\Delta (v)=(1,2,...,p-c+1,N(3,p)-s_p-q+1-N(3,p-q),1,1,...,1,0).$$ Therefore, in order to prove that $\Delta (v)$ is an $O$-sequence, it remains to show that \begin{equation}\label{qqq}p-q+2\geq N(3,p)-s_p-q+1-N(3,p-q)\geq 1.\end{equation}
The second inequality of (\ref{qqq}) follows immediately from the definition of $q$. In order to prove the first inequality, notice that, by definition of $q$, we have $$N(3,p)-s_p-(q-1)<N(3,p-(q-1)).$$ Hence $$N(3,p)-s_p<q-1+N(3,p-q+1)$$$$=q-1+N(3,p-q)+p-q+2=N(3,p-q)+p+1,$$ and the first inequality of (\ref{qqq}) follows.\\
This completes the proof of the claim and that of the theorem.{\ }{\ }\qed \\
\\
\\
{\bf \huge References}\\
\\
$[BG]$ {\ } A.M. Bigatti and A.V. Geramita: {\it Level Algebras, Lex Segments and Minimal Hilbert Functions}, Comm. in Algebra, 31, 1427-1451, 2003.\\
$[CI]$ {\ } Y.H. Cho and A. Iarrobino: {\it Inverse Systems of Zero-dimensional Schemes in ${\bf P}^n$}, preprint.\\
$[EI]$ {\ } J. Emsalem and A. Iarrobino: {\it Some zero-dimensional generic singularities; finite algebras having small tangent space}, Compositio Math., 36, 145-188, 1978.\\
$[FL]$ {\ } R. Fr\"oberg and D. Laksov: {\it Compressed Algebras}, Conference on Complete Intersections in Acireale, Lecture Notes in Mathematics, 1092, 121-151, Springer-Verlag, 1984.\\
$[Ge]$ {\ } A.V. Geramita: {\it Inverse Systems of Fat Points: Waring's Problem, Secant Varieties and Veronese Varieties and Parametric Spaces of Gorenstein Ideals}, Queen's Papers in Pure and Applied Mathematics, No. 102, The Curves Seminar at Queen's, Vol. X, 3-114, 1996.\\
$[Ia]$ {\ } A. Iarrobino: {\it Compressed Algebras: Artin algebras having given socle degrees and maximal length}, Trans. Amer. Math. Soc., 285, 337-378, 1984.\\
$[IK]$ {\ } A. Iarrobino, V. Kanev: {\it Power Sums, Gorenstein Algebras, and Determinantal Loci}, Springer Lecture Notes in Mathematics, No. 1721, Springer, Heidelberg, 1999.\\
$[St]$ {\ } R. Stanley: {\it Hilbert functions of graded algebras}, Adv. Math., 28, 57-83, 1978.\\
$[Za]$ {\ } F. Zanello: {\it Extending the idea of compressed algebra to arbitrary socle-vectors}, accepted for publication on J. of Algebra.

}

\end{document}